\documentstyle[11pt]{article}
\title{Gromov-Witten Invariants in Algebraic Geometry}
\author{K. Behrend}
\catcode`\@=11

\font\tenmsx=msam10
\font\sevenmsx=msam7
\font\fivemsx=msam5
\font\tenmsy=msbm10
\font\sevenmsy=msbm7
\font\fivemsy=msbm5
\newfam\msxfam
\newfam\msyfam
\textfont\msxfam=\tenmsx  \scriptfont\msxfam=\sevenmsx
  \scriptscriptfont\msxfam=\fivemsx
\textfont\msyfam=\tenmsy  \scriptfont\msyfam=\sevenmsy
  \scriptscriptfont\msyfam=\fivemsy

\def\hexnumber@#1{\ifcase#1 0\or1\or2\or3\or4\or5\or6\or7\or8\or9\or
        A\or B\or C\or D\or E\or F\fi }

\font\teneuf=eufm10
\font\seveneuf=eufm7
\font\fiveeuf=eufm5
\newfam\euffam
\textfont\euffam=\teneuf
\scriptfont\euffam=\seveneuf
\scriptscriptfont\euffam=\fiveeuf
\def\frak{\ifmmode\let\next\frak@\else
 \def\next{\Err@{Use \string\frak\space only in math mode}}\fi\next}
\def\goth{\relaxnext@\ifmmode\let\next\frak@\else
 \def\next{\Err@{Use \string\goth\space only in math mode}}\fi\next}
\def\frak@#1{{\frak@@{#1}}}
\def\frak@@#1{\fam\euffam#1}

\edef\msx@{\hexnumber@\msxfam}
\edef\msy@{\hexnumber@\msyfam}

\mathchardef\boxdot="2\msx@00
\mathchardef\boxplus="2\msx@01
\mathchardef\boxtimes="2\msx@02
\mathchardef\square="0\msx@03
\mathchardef\blacksquare="0\msx@04
\mathchardef\centerdot="2\msx@05
\mathchardef\lozenge="0\msx@06
\mathchardef\blacklozenge="0\msx@07
\mathchardef\circlearrowright="3\msx@08
\mathchardef\circlearrowleft="3\msx@09
\mathchardef\rightleftharpoons="3\msx@0A
\mathchardef\leftrightharpoons="3\msx@0B
\mathchardef\boxminus="2\msx@0C
\mathchardef\Vdash="3\msx@0D
\mathchardef\Vvdash="3\msx@0E
\mathchardef\vDash="3\msx@0F
\mathchardef\twoheadrightarrow="3\msx@10
\mathchardef\twoheadleftarrow="3\msx@11
\mathchardef\leftleftarrows="3\msx@12
\mathchardef\rightrightarrows="3\msx@13
\mathchardef\upuparrows="3\msx@14
\mathchardef\downdownarrows="3\msx@15
\mathchardef\upharpoonright="3\msx@16

\mathchardef\downharpoonright="3\msx@17
\mathchardef\upharpoonleft="3\msx@18
\mathchardef\downharpoonleft="3\msx@19
\mathchardef\rightarrowtail="3\msx@1A
\mathchardef\leftarrowtail="3\msx@1B
\mathchardef\leftrightarrows="3\msx@1C
\mathchardef\rightleftarrows="3\msx@1D
\mathchardef\Lsh="3\msx@1E
\mathchardef\Rsh="3\msx@1F
\mathchardef\rightsquigarrow="3\msx@20
\mathchardef\leftrightsquigarrow="3\msx@21
\mathchardef\looparrowleft="3\msx@22
\mathchardef\looparrowright="3\msx@23
\mathchardef\circeq="3\msx@24
\mathchardef\succsim="3\msx@25
\mathchardef\gtrsim="3\msx@26
\mathchardef\gtrapprox="3\msx@27
\mathchardef\multimap="3\msx@28
\mathchardef\therefore="3\msx@29
\mathchardef\because="3\msx@2A
\mathchardef\doteqdot="3\msx@2B

\mathchardef\triangleq="3\msx@2C
\mathchardef\precsim="3\msx@2D
\mathchardef\lesssim="3\msx@2E
\mathchardef\lessapprox="3\msx@2F
\mathchardef\eqslantless="3\msx@30
\mathchardef\eqslantgtr="3\msx@31
\mathchardef\curlyeqprec="3\msx@32
\mathchardef\curlyeqsucc="3\msx@33
\mathchardef\preccurlyeq="3\msx@34
\mathchardef\leqq="3\msx@35
\mathchardef\leqslant="3\msx@36
\mathchardef\lessgtr="3\msx@37
\mathchardef\backprime="0\msx@38
\mathchardef\risingdotseq="3\msx@3A
\mathchardef\fallingdotseq="3\msx@3B
\mathchardef\succcurlyeq="3\msx@3C
\mathchardef\geqq="3\msx@3D
\mathchardef\geqslant="3\msx@3E
\mathchardef\gtrless="3\msx@3F
\mathchardef\sqsubset="3\msx@40
\mathchardef\sqsupset="3\msx@41
\mathchardef\vartriangleright="3\msx@42
\mathchardef\vartriangleleft="3\msx@43
\mathchardef\trianglerighteq="3\msx@44
\mathchardef\trianglelefteq="3\msx@45
\mathchardef\bigstar="0\msx@46
\mathchardef\between="3\msx@47
\mathchardef\blacktriangledown="0\msx@48
\mathchardef\blacktriangleright="3\msx@49
\mathchardef\blacktriangleleft="3\msx@4A
\mathchardef\vartriangle="0\msx@4D
\mathchardef\blacktriangle="0\msx@4E
\mathchardef\triangledown="0\msx@4F
\mathchardef\eqcirc="3\msx@50
\mathchardef\lesseqgtr="3\msx@51
\mathchardef\gtreqless="3\msx@52
\mathchardef\lesseqqgtr="3\msx@53
\mathchardef\gtreqqless="3\msx@54
\mathchardef\Rrightarrow="3\msx@56
\mathchardef\Lleftarrow="3\msx@57
\mathchardef\veebar="2\msx@59
\mathchardef\barwedge="2\msx@5A
\mathchardef\doublebarwedge="2\msx@5B
\mathchardef\angle="0\msx@5C
\mathchardef\measuredangle="0\msx@5D
\mathchardef\sphericalangle="0\msx@5E
\mathchardef\varpropto="3\msx@5F
\mathchardef\smallsmile="3\msx@60
\mathchardef\smallfrown="3\msx@61
\mathchardef\Subset="3\msx@62
\mathchardef\Supset="3\msx@63
\mathchardef\Cup="2\msx@64

\mathchardef\Cap="2\msx@65

\mathchardef\curlywedge="2\msx@66
\mathchardef\curlyvee="2\msx@67
\mathchardef\leftthreetimes="2\msx@68
\mathchardef\rightthreetimes="2\msx@69
\mathchardef\subseteqq="3\msx@6A
\mathchardef\supseteqq="3\msx@6B
\mathchardef\bumpeq="3\msx@6C
\mathchardef\Bumpeq="3\msx@6D
\mathchardef\lll="3\msx@6E

\mathchardef\ggg="3\msx@6F

\mathchardef\circledS="0\msx@73
\mathchardef\pitchfork="3\msx@74
\mathchardef\dotplus="2\msx@75
\mathchardef\backsim="3\msx@76
\mathchardef\backsimeq="3\msx@77
\mathchardef\complement="0\msx@7B
\mathchardef\intercal="2\msx@7C
\mathchardef\circledcirc="2\msx@7D
\mathchardef\circledast="2\msx@7E
\mathchardef\circleddash="2\msx@7F
\def\ulcorner{\delimiter"4\msx@70\msx@70 }
\def\urcorner{\delimiter"5\msx@71\msx@71 }
\def\llcorner{\delimiter"4\msx@78\msx@78 }
\def\lrcorner{\delimiter"5\msx@79\msx@79 }
\def\yen{\mathhexbox\msx@55 }
\def\checkmark{\mathhexbox\msx@58 }
\def\circledR{\mathhexbox\msx@72 }
\def\maltese{\mathhexbox\msx@7A }
\mathchardef\lvertneqq="3\msy@00
\mathchardef\gvertneqq="3\msy@01
\mathchardef\nleq="3\msy@02
\mathchardef\ngeq="3\msy@03
\mathchardef\nless="3\msy@04
\mathchardef\ngtr="3\msy@05
\mathchardef\nprec="3\msy@06
\mathchardef\nsucc="3\msy@07
\mathchardef\lneqq="3\msy@08
\mathchardef\gneqq="3\msy@09
\mathchardef\nleqslant="3\msy@0A
\mathchardef\ngeqslant="3\msy@0B
\mathchardef\lneq="3\msy@0C
\mathchardef\gneq="3\msy@0D
\mathchardef\npreceq="3\msy@0E
\mathchardef\nsucceq="3\msy@0F
\mathchardef\precnsim="3\msy@10
\mathchardef\succnsim="3\msy@11
\mathchardef\lnsim="3\msy@12
\mathchardef\gnsim="3\msy@13
\mathchardef\nleqq="3\msy@14
\mathchardef\ngeqq="3\msy@15
\mathchardef\precneqq="3\msy@16
\mathchardef\succneqq="3\msy@17
\mathchardef\precnapprox="3\msy@18
\mathchardef\succnapprox="3\msy@19
\mathchardef\lnapprox="3\msy@1A
\mathchardef\gnapprox="3\msy@1B
\mathchardef\nsim="3\msy@1C
\mathchardef\ncong="3\msy@1D

\mathchardef\varsubsetneq="3\msy@20
\mathchardef\varsupsetneq="3\msy@21
\mathchardef\nsubseteqq="3\msy@22
\mathchardef\nsupseteqq="3\msy@23
\mathchardef\subsetneqq="3\msy@24
\mathchardef\supsetneqq="3\msy@25
\mathchardef\varsubsetneqq="3\msy@26
\mathchardef\varsupsetneqq="3\msy@27
\mathchardef\subsetneq="3\msy@28
\mathchardef\supsetneq="3\msy@29
\mathchardef\nsubseteq="3\msy@2A
\mathchardef\nsupseteq="3\msy@2B
\mathchardef\nparallel="3\msy@2C
\mathchardef\nmid="3\msy@2D
\mathchardef\nshortmid="3\msy@2E
\mathchardef\nshortparallel="3\msy@2F
\mathchardef\nvdash="3\msy@30
\mathchardef\nVdash="3\msy@31
\mathchardef\nvDash="3\msy@32
\mathchardef\nVDash="3\msy@33
\mathchardef\ntrianglerighteq="3\msy@34
\mathchardef\ntrianglelefteq="3\msy@35
\mathchardef\ntriangleleft="3\msy@36
\mathchardef\ntriangleright="3\msy@37
\mathchardef\nleftarrow="3\msy@38
\mathchardef\nrightarrow="3\msy@39
\mathchardef\nLeftarrow="3\msy@3A
\mathchardef\nRightarrow="3\msy@3B
\mathchardef\nLeftrightarrow="3\msy@3C
\mathchardef\nleftrightarrow="3\msy@3D
\mathchardef\divideontimes="2\msy@3E
\mathchardef\varnothing="0\msy@3F
\mathchardef\nexists="0\msy@40
\mathchardef\mho="0\msy@66
\mathchardef\eth="0\msy@67
\mathchardef\eqsim="3\msy@68
\mathchardef\beth="0\msy@69
\mathchardef\gimel="0\msy@6A
\mathchardef\daleth="0\msy@6B
\mathchardef\lessdot="3\msy@6C
\mathchardef\gtrdot="3\msy@6D
\mathchardef\ltimes="2\msy@6E
\mathchardef\rtimes="2\msy@6F
\mathchardef\shortmid="3\msy@70
\mathchardef\shortparallel="3\msy@71
\mathchardef\smallsetminus="2\msy@72
\mathchardef\thicksim="3\msy@73
\mathchardef\thickapprox="3\msy@74
\mathchardef\approxeq="3\msy@75
\mathchardef\succapprox="3\msy@76
\mathchardef\precapprox="3\msy@77
\mathchardef\curvearrowleft="3\msy@78
\mathchardef\curvearrowright="3\msy@79
\mathchardef\digamma="0\msy@7A
\mathchardef\varkappa="0\msy@7B
\mathchardef\hslash="0\msy@7D
\mathchardef\hbar="0\msy@7E
\mathchardef\backepsilon="3\msy@7F
\def\Bbb{\ifmmode\let\next\Bbb@\else
 \def\next{\errmessage{Use \string\Bbb\space only in math mode}}\fi\next}
\def\Bbb@#1{{\Bbb@@{#1}}}
\def\Bbb@@#1{\fam\msyfam#1}

\catcode`\@=12

\newtheorem{prop}{Proposition}
\newtheorem{lem}[prop]{Lemma}

\newtheorem{them}[prop]{Theorem}

\newtheorem{defnp}[prop]{Definition}
\newtheorem{numconp}[prop]{Construction}
\newtheorem{numrmkp}[prop]{Remark}
\newtheorem{numrmksp}[prop]{Remarks}
\newtheorem{conp}[prop]{}

\newtheorem{warningp}{Warning}
\newtheorem{notep}{Note}
\newtheorem{claimp}{Claim}
\newtheorem{examplep}{Example}
\newtheorem{examplesp}{Examples}
\newtheorem{rmkp}{Remark}
\newtheorem{rmksp}{Remarks}

\newenvironment{rmk}{\begin{rmkp}\rm}{\end{rmkp}}

\newenvironment{pf}{\begin{trivlist}\item[]{\sc Proof.}}%
            {\nolinebreak $\Box$ \end{trivlist}}

\newcommand{\noprint}[1]{}

\renewcommand{\tilde}{\widetilde}
\newcommand{\dual}[1]{{{#1}^{\vee}}}

\newcommand{\virt}{\mbox{\tiny virt}}

\newcommand{\upst}{^{\ast}}

\newcommand{\upsh}{^{!}}
\newcommand{\lst}{_{\ast}}

\newcommand{\com}{^{\scriptscriptstyle\bullet}}

\newcommand{\TT}{{\frak T}}

\newcommand{\PP}{{\frak P}}

\newcommand{\EE}{{\frak E}}

\newcommand{\GG}{{\frak G}}
\newcommand{\NN}{{\frak N}}
\newcommand{\MM}{{\frak M}}
\newcommand{\HH}{{\frak H}}

\newcommand{\Gg}{{\frak g}}

\newcommand{\Pp}{{\frak p}}

\newcommand{\Nn}{{\frak n}}

\newcommand{\CC}{{\frak C}}
\newcommand{\zz}{{\Bbb Z}}

\newcommand{\qq}{{\Bbb Q}}
\newcommand{\pp}{{\Bbb P}}

\newcommand{\tT}{{\cal T}}

\renewcommand{\O}{{\cal O}}

\newcommand{\cC}{{\cal C}}

\newcommand{\del}{\partial}
\newcommand{\resto}{{ \mid }}

\newcommand{\rk}{\mathop{\rm rk}}

\newcommand{\Mor}{\mathop{\rm Mor}\nolimits}

\newcommand{\multdeg}{\mathop{\rm\textstyle \mbox{mult-deg}}\nolimits}

\newcommand{\id}{\mathop{\rm id}}

\newcommand{\Pic}{\mathop{\rm Pic}\nolimits}

\newcommand{\comp}{\mathbin{{\scriptstyle\circ}}}
\newcommand{\ol}{\overline}

\newcommand{\ldiag}[1]%
       {\makebox[0cm]{${\scriptstyle#1}\downarrow\phantom{\scriptstyle#1}$}}
\newcommand{\ldiagup}[1]%
       {\makebox[0cm]{${\scriptstyle#1}\uparrow\phantom{\scriptstyle#1}$}}
\newcommand{\rdiag}[1]%
       {\makebox[0cm]{$\phantom{\scriptstyle#1}\downarrow{\scriptstyle#1}$}}
\newcommand{\sediagr}[1]%
       {\makebox[0cm]{$\phantom{\scriptstyle#1}\searrow{\scriptstyle#1}$}}
\newcommand{\nediagr}[1]%
       {\makebox[0cm]{$\phantom{\scriptstyle#1}\nearrow{\scriptstyle#1}$}}
\newcommand{\rdiagup}[1]%
       {\makebox[0cm]{$\phantom{\scriptstyle#1}\uparrow{\scriptstyle#1}$}}
\newcommand{\swdiag}[1]%
       {\makebox[0cm]{$\phantom{\scriptstyle#1}\swarrow{\scriptstyle#1}$}}
\newcommand{\sediag}[1]%
       {\makebox[0cm]{${\scriptstyle#1}\searrow\phantom{\scriptstyle#1}$}}
\newcommand{\nediag}[1]%
       {\makebox[0cm]{${\scriptstyle#1}\nearrow\phantom{\scriptstyle#1}$}}

\newcommand{\longiso}{\stackrel{\textstyle\sim}{\longrightarrow}}
\newcommand{\iso}{\stackrel{\sim}{\rightarrow}}
\newcommand{\longisoback}{\stackrel{\textstyle\sim}{\longleftarrow}}

\newcommand{\comdia}[9]{%
\begin{array}{ccc}
#1 & \stackrel{#2}{\longrightarrow} & #3 \\
\ldiag{#4} & #5 & \rdiag{#6} \\
#7 & \stackrel{#8}{\longrightarrow} & #9
\end{array}}

\newcommand{\comtri}[6]{%
\begin{array}{ccc}
#1 & \stackrel{#2}{\longrightarrow} & #3         \\
   & \sediag{#4}                    & \rdiag{#5} \\
   &                                & #6
\end{array}}

\newcommand{\inversecomtri}[6]{%
\begin{array}{ccc}
#1           &                                &            \\
\ldiag{#2}   & \sediagr{#3}                   &            \\
#4           & \stackrel{#5}{\longrightarrow} & #6
\end{array}}

\newcommand{\overtoparrow}%
{\makebox[0cm]{\beginpicture
\setcoordinatesystem units <.8cm,.4cm> point at 0 0
\setplotarea x from -3 to 3, y from 0 to 1
\setquadratic
\plot -3 0 0 1 3 0 /
\put{\vector(3,-1){0}}[Bl] at 3 0
\endpicture}}

\newcommand{\underbottomarrow}%
{\makebox[0cm]{\beginpicture
\setcoordinatesystem units <.8cm,.4cm> point at 0 0
\setplotarea x from -3 to 3, y from 0 to 1
\setquadratic
\plot -3 1 0 0 3 1 /
\put{\vector(3,1){0}}[Bl] at 3 1
\endpicture}}

\renewcommand{\t}{_{\tau}}

\newcommand{\ses}[5]%
{0\longrightarrow#1\stackrel{#2}{ \longrightarrow}#3\stackrel{#4}{
\longrightarrow}#5\longrightarrow0}

\newcommand{\dt}[6]%
{#1\stackrel{#2}{\longrightarrow}#3 \stackrel{#4}{\longrightarrow}#5
\stackrel{#6}{\longrightarrow} #1[1]}

\begin{document} \sloppy
\date{January 13, 1996}
\maketitle
\begin{abstract}
{Gromov-Witten invariants for arbitrary projective varieties
and arbitrary genus are constructed using the techniques from [K.
Behrend, B. Fantechi. {\em The Intrinsic Normal Cone.\/}]}
\end{abstract}

\subsection{Introduction}

In \cite{BM} the problem of constructing the Gromov-Witten invariants
of a smooth projective variety $V$ was reduced to defining a `virtual
fundamental class'
$$[\ol{M}_{g,n}(V,\beta)]^{\virt}\in A_{(1-g)(\dim
V-3)-\beta(\omega_V)+n}(\ol{M}_{g,n}(V,\beta))$$
in the Chow group of the algebraic stack
$$\ol{M}_{g,n}(V,\beta)$$
of stable maps of class $\beta\in H_2(V)$ from an $n$-marked prestable
curve of genus $g$ to $V$.

If $g=0$ and $V$ is convex (i.e.\ $H^1(\pp^1,f\upst T_V)=0$, for all
$f:\pp^1\to V$), then $\ol{M}_{0,n}(V,\beta)$ is smooth of the
expected dimension $\dim V-3-\beta(\omega_V)+n$ and the usual
fundamental class
$$[\ol{M}_{g,n}(V,\beta)]$$
will work. This was proved in \cite{BM}.

In this paper we treat the general case using the construction from
\cite{BF}. Recall from [ibid.] that virtual fundamental classes are
constructed using an {\em obstruction theory}, and the {\em intrinsic
normal cone}. The obstruction theory serves to give rise to a vector
bundle stack $\EE$, into which the intrinsic normal cone $\CC$ can be
embedded as a closed subcone stack. The virtual fundamental class is
then obtained by intersecting $\CC$ with the zero section of $\EE$.

In our context, this process works as follows. Let $\MM_{g,n}$ be the
algebraic stack of $n$-marked prestable curves of genus $g$. This is
an algebraic stack, not of Deligne-Mumford (or even finite) type, but
smooth of dimension $3(g-1)+n$. There is a canonical morphism
$$\ol{M}_{g,n}(V,\beta) \to \MM_{g,n},$$ given by forgetting the map,
retaining the curve (but not stabilizing). Then $\ol{M}_{g,n}(V,\beta)
\to \MM_{g,n}$ is an open substack of a stack of morphisms, and as
such has a relative obstruction theory, which in this case is
$\dual{(R\pi\lst f\upst T_V)}$, where $\pi:C\to\ol{M}_{g,n}(V,\beta)$
is the universal curve and $f:C\to V$ is the universal stable
map. Saying that $\dual{(R\pi\lst f\upst T_V)}$is a relative
obstruction theory means that there is a homomorphism
$$\phi:\dual{(R\pi\lst f\upst T_V)} \longrightarrow
L\com_{\ol{M}_{g,n}(V,\beta) / \MM_{g,n}},$$ (where $L\com$ is the
cotangent complex) such that $h^0(\phi)$ is an isomorphism and
$h^{-1}(\phi)$ is surjective.

The homomorphism $\phi$ induces a closed immersion
\[\dual{\phi}:\NN_{\ol{M}_{g,n}(V,\beta) / \MM_{g,n}}\longrightarrow
h^1/h^0(R\pi\lst f\upst T_V)$$ of abelian cone stacks (see \cite{BF})
over $\ol{M}_{g,n}(V,\beta)$, where $\NN$ is the relative intrinsic
normal sheaf. The relative intrinsic normal cone
$\CC_{\ol{M}_{g,n}(V,\beta) / \MM_{g,n}}$ is a closed subcone stack of
$\NN_{\ol{M}_{g,n}(V,\beta) / \MM_{g,n}}$, and so we get a closed
immersion of cone stacks $$\CC_{\ol{M}_{g,n}(V,\beta) /
\MM_{g,n}} \longrightarrow h^1/h^0(R\pi\lst f\upst T_V).$$ Now since
$R\pi\lst f\upst T_V$ has global resolutions (see
Proposition~\ref{egr}), we may intersect $\CC_{\ol{M}_{g,n}(V,\beta) /
\MM_{g,n}}$ with the zero section of the vector bundle stack
$h^1/h^0(R\pi\lst f\upst T_V)$ to get the virtual fundamental class
$[\ol{M}_{g,n}(V,\beta)]^{\virt}$.

The fundamental axioms (see \cite{KM}) Gromov-Witten invariants need
to satisfy to deserve their name are reduced in \cite{BM} to five
basic compatibilities between the virtual fundamental classes. These
follow from the basic properties proved in \cite{BF}. The dimension
axiom, for example, follows from the basic fact that the intrinsic
normal cone always has dimension zero.

We also show that if $V=G/P$, for a reductive group $G$ and a parabolic
subgroup $P$, there is an alternative construction of the virtual
fundamental classes avoiding the intrinsic normal cone. We construct a
cone $C$ in the vector bundle $R^1\pi\lst\O\otimes\Gg$ on
$\ol{M}_{g,n}(V,\beta)$, which may then be intersected with the zero
section of $R^1\pi\lst\O\otimes\Gg$ to obtain the virtual fundamental
class. This cone $C$ is constructed as the normal cone of an embedding
of $\ol{M}_{g,n}(V,\beta)$ into a certain stack of principal
$P$-bundles (which is smooth, but not of Deligne-Mumford type).

\smallskip

A construction of Gromov-Witten invariants using a cone inside a vector
bundle has also been announced by J. Li and G. Tian. Their methods
differ from ours in that they use analytic methods, including the
Kuranishi map.

\smallskip

Most of this work was done during a stay at the Max-Planck-Institut
f\"ur Mathematik in Bonn, and I would like to take this opportunity to
acknowledge the hospitality and the wonderful atmosphere at the MPI\@.
I would also like to thank Professors G. Harder and Yu.\ Manin for fruitful
discussions about the constructions in this paper.

\subsection{Preliminaries on Prestable Curves}

Let $k$ be a field. We shall work over the category of locally noetherian
$k$-schemes (with the fppf-topology). For a modular graph $\tau$ (see
\cite{BM},
Definition~1.5) let $\MM(\tau)$ denote the $k$-stack of $\tau$-marked prestable
curves (which are defined in \cite{BM}, Definition~2.6).

\begin{lem}\label{soss}
The algebraic $k$-stack $\ol{M}(\tau)$ of stable $\tau$-marked curves
is an open substack of $\MM(\tau)$.
\end{lem}
\begin{pf}
Let $\cC_v\rightarrow\MM(\tau)$ be the universal curve
corresponding to the vertex $v\in V_{\tau}$. Let $\tilde{\cC}_v$ be
the stabilization. Then $\ol{M}(\tau)$ is the substack of $\MM(\tau)$
over which all $p_v:\cC_v\to\tilde{\cC}_v$ are isomorphisms. This is
open because the $\cC_v$ are proper over $\MM(\tau)$.
\end{pf}

Now consider a modular graph $\tau'$ obtained from $\tau$ by adding
some tails. We get an induced morphism of $k$-stacks
$\MM(\tau')\to\MM(\tau)$ which simply forgets the markings
corresponding to the tails $S_{\tau'}-S_{\tau}$. If
$S_{\tau'}-S_{\tau}$ has cardinality 1, then
$\MM(\tau')\rightarrow\MM(\tau)$ is a smooth
curve, hence representable and smooth of relative dimension 1. So by
induction, $\MM(\tau')\to\MM(\tau)$ is representable and smooth of
relative dimension $\#(S_{\tau'}-S_{\tau})$. By Lemma~\ref{soss} the
same is true for $\ol{M}(\tau')\to\MM(\tau)$.

\begin{prop}  The stack $\MM(\tau)$ is a smooth algebraic $k$-stack of
dimension
\[\dim(\tau)=\#S_{\tau}-\#E_{\tau}-3\chi(\tau).\]
\end{prop}
\begin{pf} For the definition of $\dim(\tau)$ and $\chi(\tau)$ see \cite{BM},
Definitions~6.1 and~6.2.

Note that for every point of $\MM(\tau)$ there exists a $\tau'$ as
above such that the induced morphism $\ol{M}(\tau')\to\MM(\tau)$
contains this given point in its image. Thus
$\coprod_{\tau'}\ol{M}(\tau')$ is a presentation of $\MM(\tau)$
showing that $\MM(\tau)$ is algebraic.
\end{pf}

Now let $\tau^s$ be the stabilization of $\tau$. Stabilization defines
a morphism of algebraic $k$-stacks
$$s:\MM(\tau)\longrightarrow\ol{M}(\tau^s).$$
If $\tau'$ is obtained as above by adjoining tails to $\tau$ such that
$\tau'$ is stable, we have a commutative diagram
\[\inversecomtri{\ol{M}(\tau')}{
}{\ol{M}(\phi)}{\MM(\tau)}{s}{\ol{M}(\tau^s).}\]
Here $\phi:\tau'\to\tau^s$ is the canonical morphism of stable modular
graphs. In fact, one may define $s$ locally by using such diagrams.
\begin{prop}
The morphism $s:\MM(\tau)\rightarrow \ol{M}(\tau^s)$ is flat.
\end{prop}
\begin{pf}
This follows by descent since the morphism $\ol{M}(\phi)$ for various
$\phi:\tau'\rightarrow \tau^s$ as above are flat.
\end{pf}

\subsection{The Virtual Fundamental Classes}

Over $\MM(\tau)$ there is a family $(\cC_v)_{v\in V_{\tau}}$ of universal
curves, with sections
$x_i:\MM(\tau)\rightarrow\cC_{\del\t(i)}$. Let $\cC(\tau)\rightarrow\MM(\tau)$
be the curve obtained from $\coprod_{v\in V\t}\cC_v$ by identifying $x_i$ and
$x_j$, for every edge $\{i,j\}\in E\t$. The curve $\cC(\tau)$ has markings
$x_i:\MM(\tau)\rightarrow\cC(\tau)$, for each $i\in S\t$. In fact,
$\cC(\tau)$ is a $\tilde{\tau}$-marked prestable curve, where $\tilde{\tau}$ is
the graph obtained from $\tau$ by contracting all edges of $\tau$. Let us
denote
the structure morphism by
\[\pi:\cC(\tau)\longrightarrow\MM(\tau).\] We shall also denote any base change
of $\pi$ by $\pi$.

Now let $V$ be a smooth projective $k$-variety, $(\tau,\beta)$ a stable
$V$-graph and let
$\Mor_{\MM(\tau)}(\tau,V)$ be the $\MM(\tau)$-space of morphisms
{}from $\cC(\tau)$ to
$V$. Denote the universal morphism by
\[f:\cC(\tau)\times\Mor_{\MM(\tau)}(\tau,V)\longrightarrow V.\] By
\cite{fgaIV} the stack $\Mor_{\MM(\tau)}(\tau,V)$ is an algebraic
$k$-stack and the
structure morphism
\[\Mor_{\MM(\tau)}(\tau,V)\longrightarrow \MM(\tau)\] is
representable.

\begin{prop}
The proper Deligne-Mumford stack $\ol{M}(V,\tau,\beta)$ of stable maps
is an open substack of $\Mor_{\MM(\tau)}(\tau,V)$.
\end{prop}
\begin{pf}
The set of points where stabilization is an isomorphism is open.
\end{pf}

To define the virtual fundamental class on $\ol{M}(V,\tau,\beta)$ we
consider the morphism $\ol{M}(V,\tau,\beta)\rightarrow\MM(\tau)$ and
denote the relative intrinsic normal cone (see \cite{BF}) by
$$\CC(V,\tau,\beta)=\CC_{\ol{M}(V,\tau,\beta)/\MM(\tau)}$$ The
intrinsic normal sheaf [ibid.] of $\ol{M}(V,\tau,\beta)$ over
$\MM(\tau)$ we shall denote by $\NN(V,\tau,\beta)$.

By the relative version of \cite{BF} Proposition~6.2 we have a perfect
relative obstruction theory [ibid.]
\[\dual{\pi\lst(\dual{e})}:\dual{R\pi\lst(f\upst T_V)}\longrightarrow
L_{\Mor_{\MM(\tau)}(\tau,V)/\MM(\tau)}\com.\] Restricting to the open
substack $\ol{M}(V,\tau,\beta)$ we get a perfect relative obstruction theory
\[\dual{\pi\lst(\dual{e})}:\dual{R\pi\lst(f\upst T_V)}\longrightarrow
L_{\ol{M}(V,\tau,\beta)/\MM(\tau)}\com,\]
which we shall also denote by $E\com(V,\tau,\beta)$.
Thus $\CC(V,\tau,\beta)$ is embedded as a closed subcone stack in the
vector bundle stack $$\EE(V,\tau,\beta)=h^1/h^0(R\pi\lst f\upst
T_V).$$

Note that the relative virtual
dimension of $\ol{M}(V,\tau,\beta)$ over $\MM(\tau)$ with respect to the
obstruction theory $\dual{R\pi\lst(f\upst T_V)}$ is equal to
\begin{eqnarray*}
\rk\dual{R\pi\lst(f\upst T_V)} & = & \chi(f\upst T_V) \\
 & = & \deg f\upst T_V+\dim V\cdot\chi(\cC(\tau)) \\
 & = & \chi(\tau)\dim V-\beta(\tau)(\omega_V).
\end{eqnarray*} Essential is the following result.

\begin{prop} \label{egr}
Let $(C,x,f)$ be a stable map over $T$ to $V$, where $T$ is a finite
type algebraic $k$-stack. Let $E$ be a vector bundle on $C$. Then
$R\pi\lst E$ has global resolutions, where $\pi:C\rightarrow T$ is the
structure map.
\end{prop}
\begin{pf}
Let $M$ be an ample invertible sheaf on $V$ and let
\[L=\omega_{C/T}(x_1+\ldots+x_n)\otimes f\upst M^{\otimes 3}.\]
By Proposition~3.9 of \cite{BM} the sheaf $L$ is ample on the fibers
of $\pi$. So for sufficiently large $N$ we have that
\begin{enumerate}
\item $\pi\upst\pi\lst(E\otimes L^{\otimes N})\rightarrow E\otimes
L^{\otimes N} $ is surjective,
\item $R^1\pi\lst(E\otimes L^{\otimes N})=0$,
\item for all $t\in T$ we have that $H^0(C_t,L^{\otimes -N}_t)=0$.
\end{enumerate}
Let $$F=\pi\upst\pi\lst(E\otimes L^{\otimes N})\otimes L^{\otimes
-N}$$ and let $H$ be the kernel of the map $F\to E$. Thus we have a
short exact sequence
\[\ses{H}{}{F}{}{E}\] of vector bundles on $C$. Note that for every
$t\in T$ we have
\begin{eqnarray*}
H^0(C_t,F) & = & H^0(C_t,\pi\lst(E\otimes L^{\otimes N})_t\otimes L_t^{\otimes
-N})  \\    & = & H^0(C_t,L_t^{\otimes -N})\otimes\pi\lst(E\otimes
L^{\otimes N})_t \\ & = & 0 \end{eqnarray*} and hence $H^0(C_t,H)=0$,
also. Therefore, $\pi\lst H$ and $\pi\lst F$ are zero and $R^1\pi\lst
H$ and $R^1\pi\lst F$ are locally free. This implies that $$R\pi\lst
E\cong [R^1\pi\lst H\rightarrow R^1\pi\lst F].$$
\end{pf}

As shown in \cite{BF}, by Proposition~\ref{egr} the obstruction theory
$\dual{R\pi\lst(f\upst T_V)}$ gives rise to a virtual fundamental
class
\[[\ol{M}(V,\tau,\beta),\dual{R\pi\lst(f\upst T_V)}]\in
A_{\dim(V,\tau,\beta)}(\ol{M}(V,\tau,\beta)),\] since
\begin{eqnarray*}
\lefteqn{\dim\MM(\tau) + \rk \dual{R\pi\lst(f\upst T_V)} }\\ & = &
\chi(\tau)(\dim V-3)-\beta(\tau)(\omega_V)+\#S\t-\#E\t \\
 & = &  \dim(V,\tau,\beta).
\end{eqnarray*}
(See Definition~6.2 in \cite{BM} for the definition of
$\dim(V,\tau,\beta)$.)

\begin{them} \label{ot} The system of virtual fundamental classes
\[J(V,\tau,\beta)=[\ol{M}(V,\tau,\beta),\dual{R\pi\lst(f\upst T_V)}]\] is an
orientation of $\ol{M}$ over $\GG_s(V)$. If $V$ is convex, on the tree level
sub-category
$\TT_s(V)$, we get back the orientation of \cite{BM}, Theorem~7.5.
\end{them}
\begin{pf} If $V$ is convex and $\tau$ a forest, then $R^1\pi\lst(f\upst
T_V)=0$, so that the virtual fundamental class is the usual
fundamental class by \cite{BF} Proposition~7.3. Thus the virtual
fundamental class agrees with the orientation of \cite{BM},
Theorem~7.5. To check that $J$ is an orientation, we need to check the
five axioms listed in
\cite{BM}, Definition~7.1. This shall be done in the next Section.
\end{pf}

\begin{rmk} As shown in \cite{BM}, we get an associated system of
Gromov-Witten
classes for $V$.
\end{rmk}

\subsection{Checking the Axioms} \label{check}

\subsubsection{{\sc Axiom I.} Mapping to a point}

Let $\tau$ be a stable $V$-graph of class zero such that $|\tau|$ is
non-empty and connected. As noted in \cite{BM} Section~7 we have
\[\ol{M}(V,\tau,0)=V\times\ol{M}(\tau)\]
which is obviously smooth over $\MM(\tau)$. In fact, the morphism
$\ol{M}(V,\tau,0)\to\MM(\tau)$ is just the composition
\[V\times\ol{M}(\tau)\longrightarrow\ol{M}(\tau)\longrightarrow\MM(\tau)\]
of projection followed by inclusion. If
$\tilde{\pi}:\cC(\tau)\to\ol{M}{\tau}$ is the universal curve over
$\ol{M}(\tau)$, then $\cC(V,\tau,0)=V\times\cC(\tau)$ and
$\pi:\cC(V,\tau,0)\to\ol{M}(V,\tau,0)$ is identified with
$\id\times\tilde{\pi}:V\times\cC(\tau)\to V\times\ol{M}(\tau)$. Hence
\begin{eqnarray*}
R^1\pi\lst f\upst T_V&=&T_V\boxtimes
R^1\tilde{\pi}\lst\O_{\cC(\tau)}\\&=&\tT^{(1)}\end{eqnarray*}
is locally free. So by \cite{BF} Proposition~7.3 we have
\begin{eqnarray*}
J(V,\tau,0) & = & c_{\rk R^1\pi\lst f\upst T_V}(R^1\pi\lst f\upst
T_V)\cdot[\ol{M}(V,\tau,0)]\\
& = & c_{g(\tau)\dim V}(\tT^{(1)})\cdot[\ol{M}(V,\tau,0)],
\end{eqnarray*}
which is Axiom~I.

\subsubsection{{\sc Axiom II.} Products}

Let $(\sigma,\alpha)$ and $(\tau,\beta)$ be stable $V$-graphs and
denote the `product' by $(\sigma\times\tau,\alpha\times\beta)$.
Note that
\[E\com(V,\sigma\times\tau,\alpha\times\beta) =
E\com(V,\sigma,\alpha)\boxplus E\com(V,\tau,\beta),\]
so by \cite{BF} Proposition~7.4 we have
\begin{eqnarray*}
J(V,\sigma\times\tau,\alpha\times\beta) & = &
[\ol{M}(V,\sigma\times\tau,\alpha\times\beta),
E\com(V,\sigma,\alpha)\boxplus E\com(V,\tau,\beta)]  \\
& = & [\ol{M}(V,\sigma,\alpha),E\com(V,\sigma,\alpha)] \times
[\ol{M}(V,\tau,\beta),E\com(V,\tau,\beta)] \\
& = & J(V,\sigma,\alpha)\times J(V,\tau,\beta),
\end{eqnarray*}
which is the product axiom.

\subsubsection{{\sc Axiom III.} Cutting Edges}

Use notation as in \cite{BM}, Section~7, modified as necessary to
avoid confusion. Let $\beta$ denote the $H_2(V)^+$-structure on both
$\sigma$ and $\tau$. Write $\MM=\MM(\tau)=\MM(\sigma)$. Consider the
cartesian diagram
\[\comdia{\ol{M}(V,\sigma,\beta)}{\ol{M}(
\Phi)}{\ol{M}(V,\tau,\beta)}{ g}{}{ }{\MM\times V}{ \Delta}{ \MM\times
V\times V}\]
of stacks over $\MM$. Let us show that the obstruction theories
$E\com(V,\tau,\beta)$ and $E\com(V,\sigma,\beta)$ are compatible over
$\Delta$ (see \cite{BF}).

Over $\ol{M}(V,\sigma,\beta)$ let us consider the following two
curves. First the curve $\cC=\cC(V,\sigma,\beta)$ obtained from the
universal curves $(C_v)_{v\in V_{\sigma}}$ by gluing according to the
edges of $\sigma$. Secondly, we have the curve $\cC'$, which we obtain
from $(C_v)_{v\in V_{\sigma}}$ by gluing according to the
edges of $\tau$. In other words,
$\cC'=\ol{M}(\Phi)\upst\cC(V,\tau,\beta)$. Moreover, $\cC$ is obtained
from $\cC'$ by identifying the two sections $x_1$ and $x_2$ of $\cC'$,
corresponding to the edge $\{i_1,i_2\}$ of $\sigma$ which is cut by
$\Phi$. Thus there is a structure morphism $p:\cC'\to\cC$ fitting into
the commutative diagram
\[\comtri{\cC'}{p}{\cC}{\pi'}{\pi}{\ol{M}(V,\sigma,\beta).}\]
We shall also use the diagram
\[\comtri{\cC'}{p}{\cC}{f'}{f}{V,}\]
where $f:\cC\to V$ is the universal map. Let $x=p\comp x_1=p\comp
x_2$.

If $E$ is any locally free sheaf on $\cC$, then for $i=1,2$ we have
the evaluation homomorphism
\[u_i:p\upst E\longrightarrow {x_i}\lst{x_i}\upst p\upst E = {x_i}\lst
x\upst E.\]
Applying $p\lst$ we get
\[p\lst(u_i):p\lst p\upst E\longrightarrow x\lst x\upst E.\]
Letting $u=p\lst(u_2)-p\lst(u_1)$ we have a short exact sequence
\[\ses{E}{}{p\lst p\upst E}{u}{x\lst x\upst E}\]
of coherent sheaves on $\cC$. Applying $R\pi\lst$ we get a
distinguished triangle
\[\dt{R\pi\lst E}{}{ R\pi'\lst p\upst E}{R\pi\lst(u)}{x\upst E}{}\]
in $D(\O_{\ol{M}(V,\sigma,\beta)})$. Taking $E=f\upst T_V$ we get the
distinguished triangle
\[\dt{R\pi\lst f\upst T_V}{}{ R\pi'\lst {f'}\upst
T_V}{R\pi\lst(u)}{x\upst f\upst T_V}{},\]
or dually,
\begin{equation}\label{esita}
\dt{x\upst f\upst \Omega_V}{\dual{R\pi\lst(u)}}{\dual{(R\pi'\lst
{f'}\upst T_V)}}{}{\dual{(R\pi\lst f\upst T_V)}}{}.
\end{equation}
Note that we have $E\com(V,\sigma,\beta)=\dual{(R\pi\lst f\upst T_V)}$
and $\ol{M}(\Phi)\upst(E\com(V,\tau,\beta))=\dual{(R\pi'\lst{f'}\upst
T_V)}$. Moreover, $L_{\Delta}\com=\Omega_V[1]\resto\MM\times V$, so
that $g\upst L_{\Delta}=x\upst f\upst \Omega_V[1]$, since $f\comp
x=p_V\comp g$. So (\ref{esita}) gives the distinguished triangle
\[{g\upst
L_{\Delta}[-1]}\stackrel{\dual{R\pi\lst(u)}}{\longrightarrow}
{\ol{M}(\Phi)\upst
E\com(V,\tau,\beta)}{\longrightarrow}{E\com(V,\sigma,\beta)}{\longrightarrow}
g\upst L_{\Delta},\]
which we may shuffle around to give
\[\dt{\ol{M}(\Phi)\upst
E\com(V,\tau,\beta)}{}{E\com(V,\sigma,\beta)}{}{g\upst
L_{\Delta}}{\dual{R\pi\lst(-u)}}.\]
Now we have the obstruction morphisms $E\com(V,\tau,\beta)\to
L\com_{\ol{M}(V,\tau,\beta)/\MM}$ and $E\com(V,\sigma,\beta)\to
L\com_{\ol{M}(V,\sigma,\beta)/\MM}$. Moreover, we have the natural
homomorphism $g\upst L_{\Delta}\to L\com_{\ol{M}(\Phi)}$. These give
rise to a homomorphism of distinguished triangles
\[\begin{array}{ccccccc}
\ol{M}(\Phi)\upst E\com(V,\tau,\beta) & \longrightarrow &
E\com(V,\sigma,\beta) & \longrightarrow & g\upst L_{\Delta} &
\stackrel{\dual{R\pi\lst(-u)}}{\longrightarrow} & \ol{M}(\Phi)\upst
E\com(V,\tau,\beta)[1]\\
\downarrow & & \downarrow & &\downarrow & &\downarrow \\
\ol{M}(\Phi)\upst L\com_{\ol{M}(V,\tau,\beta)/\MM} & \longrightarrow &
L\com_{\ol{M}(V,\sigma,\beta)/\MM} & \longrightarrow &
L\com_{\ol{M}(\Phi)} & \longrightarrow & \ol{M}(\Phi)\upst
L\com_{\ol{M}(V,\tau,\beta)/\MM}[1],
\end{array}\]
showing that $E\com(V,\tau,\beta)$ and $E\com(V,\sigma,\beta)$ are
compatible over $\Delta$.
Hence by \cite{BF} Proposition~7.5 we have
\[\Delta\upsh J(V,\tau,\beta)=J(V,\sigma,\beta)\]
which is Axiom~III.

\subsubsection{{\sc Axiom IV.} Forgetting Tails}

Let us deal with the incomplete case, leaving the tripod losing cases
to the reader.
Letting $\cC\to\MM(\tau)$ be the universal curve corresponding to the
vertex $w\in V_{\tau}$ (notation from \cite{BM}, Section~7). We have a
cartesian diagram of algebraic $k$-stacks
\[\comdia{\ol{M}(V,\sigma,\beta)}{\ol{M}(\Phi)}{\ol{M}(V,\tau,\beta)}{
d}{ }{ }{\cC}{}{\MM(\tau).}\]
By \cite{BF} Proposition~7.2 we have
\[\ol{M}(\Phi)\upst
J(V,\tau,\beta) = [\ol{M}(V,\sigma,\beta),\ol{M}(\Phi)\upst
E\com(V,\tau,\beta)].\]
Here the class on the right hand side is the virtual fundamental class
defined by the relative intrinsic normal cone of the morphism $d$ and
the relative obstruction theory $\ol{M}(\Phi)\upst
E\com(V,\tau,\beta)$.
Note that the structure morphism
$\ol{M}(V,\sigma,\beta)\to\MM(\sigma)$ factors through
$d:\ol{M}(V,\sigma,\beta)\to\cC$.
\[\comtri{\ol{M}(V,\sigma,\beta)}{d}{\cC}{}{}{\MM(\sigma)}\]

The morphism $d:\ol{M}(V,\sigma,\beta)\to\cC$ associates to the stable
map $(C,x,h)$ the pair $((C',x'),y)$, where $(C',x',h')$ is the image
of $(C,x,h)$ under $\ol{M}(\Phi)$ and $(C',x')$ the underlying
$\tau$-marked prestable curve. Letting $x_f$ be the section of $C_v$
corresponding to the flag $f$, we obtain $(C',x',h')$ by forgetting
$x_f$ and stabilizing. Moreover, $y$ is the image of the forgotten
section $x_f$ in $C_w'$.

The morphism $\cC\to\MM(\sigma)$ associates to the pair $((C,x),y)$,
where $(C,x)$ is a $\tau$-marked prestable curve and $y$ a section of
$C_w$, the $\sigma$-marked prestable curve $(\tilde{C},\tilde{x})$
obtained as follows. For $v'\not=v$ we have $\tilde{C}_{v'}=C_{w'}$,
where $w'$ is the vertex of $\tau$ corresponding to $v'$. The curve
$(\tilde{C}_v,(\tilde{x}_j)_{j\in F_{\sigma}(v)})$ is obtained from
$((C_w,(x_j)_{j\in F_{\tau}(w)}),y)$ by `prestabilizing' (i.e.\
separating the special points) as in \cite{knudsen}, Definition~2.3.

\begin{lem}\label{loem}
The morphism $\cC\to\MM(\sigma)$ is \'etale.
\end{lem}
\begin{pf}
We will use the formal criterion for \'etaleness. Without loss of
generality assume that $w$ is the only vertex of $\tau$. So let
$((C,x),y)$ be a $\tau$-marked prestable curve with section over the
scheme $T$, $T\to T'$ a square zero extension and $(C',x')$ a
$\sigma$-marked prestable curve over $T'$ such that $(C',x')\resto T$
is the prestabilization of $((C,x),y)$. We may assume that we may
choose additional sections $s$ of $C$ over $T$, making $(C,x,s)$ a
stable marked curve. Then we extend the sections $s$ to sections $s'$
of $C'$ over $T'$. Taking the stabilization of $(C',x',s')$ after
forgetting the section $x_f'$ gives an extension of $((C,x),y)$ to
$T'$ whose prestabilization is $(C',x')$.
\end{pf}

Consider the natural morphism
$p:\cC(V,\sigma,\beta)\rightarrow\ol{M}(\Phi)\upst\cC(V,\tau)$, which
fits into the two commutative diagrams
\[\comtri{\cC(V,\sigma,\beta)}{p}{\ol{M}(\Phi)\upst \cC(V,\tau,\beta)}
{\pi} {\pi'} {\ol{M}(V,\sigma,\beta)}\]
and
\[\comtri{\cC(V,\sigma,\beta)}{p}{\ol{M}(\Phi)\upst \cC(V,\tau,\beta)}
{f}{f'}{\quad V\quad.}\]
Whenever $E$ is a locally free sheaf on
$\ol{M}(\Phi)\upst\cC(V,\tau,\beta)$ the canonical homomorphism
$E\rightarrow p\lst p\upst E$ is an isomorphism. Applying this
principle to $E={f'}\upst T_V$ we get an isomorphism
\[{f'}\upst T_V\longrightarrow p\lst f\upst T_V.\]
Applying $R\pi'\lst$ to this, gives an isomorphism
\[R\pi'\lst {f'}\upst T_V\longrightarrow R\pi\lst f\upst T_V.\]
Noting that $R\pi'\lst{f'}\upst T_V=\ol{M}(\Phi)\upst
E\com(V,\tau,\beta)$ we get an isomorphism
\[\ol{M}(\Phi)\upst E\com(V,\tau,\beta)\longrightarrow
E\com(V,\sigma,\beta)\]
and whence an isomorphism
\[\EE(V,\sigma,\beta)\longrightarrow\ol{M}(\Phi)\upst
\EE(V,\tau,\beta).\]
By \cite{BF} Proposition~7.1 there is a natural isomorphism
\[\CC_{\ol{M}(V,\sigma,\beta)/\cC}\longrightarrow\ol{M}(\Phi)\upst
\CC_{\ol{M}(V,\tau,\beta)/\MM(\tau)}.\]
By Lemma~\ref{loem} we have a canonical isomorphism
\[\CC_{\ol{M}(V,\sigma,\beta)/\cC}\longrightarrow
\CC_{\ol{M}(V,\sigma,\beta)/\MM(\sigma)},\]
such that the diagram
\[\begin{array}{ccc}
\CC_{\ol{M}(V,\sigma,\beta)/\MM(\sigma)} & \longisoback &
\CC_{\ol{M}(V,\sigma,\beta)/\cC} \\
\cap & & \cap\\
\EE(V,\sigma,\beta) & \longiso  & \ol{M}(\Phi)\upst
\EE(V,\tau,\beta)\end{array}\]
commutes. So finally, we have
\begin{eqnarray*}
\ol{M}(\Phi)\upst J(V,\tau,\beta) & = & [\ol{M}(V,\sigma,\beta),
\ol{M}(\Phi)\upst E\com(V,\tau,\beta)]\\
& = & [\ol{M}(V,\sigma,\beta),E\com(V,\sigma,\beta)]\\
& = & J(V,\sigma,\beta),
\end{eqnarray*}
which is Axiom~IV.

\subsubsection{{\sc Axiom V.} Isogenies}

Before we start with the proof, some general remarks. Let
$\Phi:\tau\to\sigma$ be an elementary contraction of stable modular
graphs, contracting the edge $\{f,\ol{f}\}$ of $\tau$.  Let
$a:\tau\rightarrow\tau'$ and $b:\sigma\rightarrow\sigma'$ be
combinatorial morphisms of modular graphs identifying $\tau$ and
$\sigma$ as the stabilizations of $\tau'$ and $\sigma'$, respectively.
Finally, let $\Phi':\tau'\to\sigma'$ be as follows. We require
$\{a(f),a(\ol{f})\}$ to be an edge of $\tau'$ and
$\Phi':\tau'\to\sigma'$ to be the elementary contraction contracting
the edge $\{a(f),a(\ol{f})\}$. Moreover, we require $\Phi$ to be the
stabilization of $\Phi'$. To fix notation, denote the vertex onto
which $\Phi'$ contracts the edge $\{a(f),a(\ol{f})\}$ by $v_0\in
V_{\sigma'}$ and let $v_1=\del_{\tau'}(a(f))$ and
$v_2=\del_{\tau'}(a(\ol{f}))$.

In this situation we get a commutative diagram of algebraic stacks
\[\comdia{\MM(\tau')}{\MM(\Phi')}{\MM(\sigma')}{ s}{}{
s}{\ol{M}(\tau)}{\ol{M}(\Phi)}{\ol{M}(\sigma).}\]
Define $\PP$ to be the fibered product
\[\comdia{\PP}{}{\MM(\sigma')}{ }{}{
s}{\ol{M}(\tau)}{\ol{M}(\Phi)}{\ol{M}(\sigma).}\]
Consider the induced morphism $l:\MM(\tau')\to\PP$.

\begin{prop} \label{evfa}
We have $l\lst[\MM(\Phi')]=s\upst[\ol{M}(\Phi)]$.
\end{prop}
\begin{pf}
First note that $\MM(\tau')$ is irreducible, since $\MM(\tau')$ is a
product of stacks of the form $\MM_{g,n}$, which are irreducible since
the stacks $\ol{M}_{g,n}$ are. Moreover, $\MM(\tau')\to\PP$ is
surjective, so that $\PP$ is irreducible, too.

Secondly, let us remark that there exist non-empty (hence dense) open
substacks $\MM(\tau')^0\subset\MM(\tau')$ and $\PP^0\subset\PP$ such
that $l$ induces an isomorphism $l^0:\MM(\tau')^0\iso\PP^0$. In fact,
let $\MM(\tau')^0$ be the open substack of $\MM(\tau')$ characterized
by the requirement that the marked curves $C_{v_1}$ and $C_{v_2}$ be
stable. To construct $\PP^0$, let $\MM(\sigma')^0$ be the open
substack of $\MM(\sigma')$ where the marked curve $C_{v_0}$ is stable.
Then set $$\PP^0=\ol{M}(\tau)\times_{\ol{M}(\sigma)}\MM(\sigma')^0.$$

These facts imply the claim.
\end{pf}

Now let $(\Phi,m):\tau\to \sigma$ be an elementary isogeny of type
forgetting a tail. Let $f\in F_{\tau}$ be the forgotten tail. Let
$a:\tau\to\tau'$ and $b:\sigma\to\sigma'$ be as above. Finally, let
$\Phi':\tau'\to\sigma'$ be the `adjoint' of a combinatorial morphism
of graphs, such that there exists a tail map $m'$, a semigroup $A$ and
$A$-structures on $\tau'$ and $\sigma'$ making $(\Phi',m')$ the
elementary isogeny of stable $A$-graphs forgetting the tail $a(f)$.
Moreover, we require $\Phi$ to be the stabilization of $\Phi'$.

Let $\PP$ be the fibered product
\[\comdia{\PP}{}{\MM(\sigma')}{}{}{s}{\ol{M}(\tau)}{}{\ol{M}(\sigma)}\]
and $\cC$ the universal curve over $\MM(\sigma')$ corresponding to
$w\in V_{\sigma'}$, where $w$ is the vertex of the forgotten tail. (If
$w$ does not exist, i.e.\ if $\Phi'$ is complete, then
$\cC=\MM(\sigma')$.) As
in the proof of Axiom~IV we have a morphism $\cC\rightarrow\MM(\tau')$
giving rise to a commutative diagram
\[\begin{array}{ccc}
\cC & \stackrel{\pi'}{\longrightarrow} & \MM(\sigma') \\
\ldiag{} & & \\
\MM(\tau') & & \rdiag{s} \\
\ldiag{s} & & \\
\ol{M}(\tau) & \stackrel{\ol{M}(\Phi)}{\longrightarrow} &
\ol{M}(\sigma)
\end{array}\]
and hence to a morphism $l:\cC\rightarrow\PP$.

\begin{prop} \label{evfat}
We have $l\lst[\pi']= s\upst [\ol{M}(\Phi)]$.
\end{prop}
\begin{pf}
Again, $\cC$ and $\PP$ are irreducible and $l$ induces an isomorphism
$l^0:\cC^0\to\PP^0$, where $\cC^0$ is the restriction of $\cC$ to
$\MM(\sigma')^0$ and
$\PP^0=\ol{M}\times_{\ol{M}(\sigma)}\MM(\sigma')^0$. Here
$\MM(\sigma')^0\subset\MM(\sigma')$ is the open substack where $C_w$
is stable.
\end{pf}

Now let us prove Axiom~V.
According to \cite{BM}, Remark~7.2, it suffices to do this
for the case that $\Phi:\tau\to \sigma$ is an elementary isogeny,
$\#J=1$ and $(a_i,\tau_i,\Phi_i)_{i\in I}$ a pullback. So we shall use
notation as in the Definition of pullback (\cite{BM}, Definition~6.10).
We shall include the $H_2(V)^+$-structures on $\sigma'$ and $\tau_i$
($i\in I$) in the notation. They shall be denoted by $\beta'$ and
$\beta_i$ ($i\in I$), respectively.
The underlying graph of $(\tau_i,\beta_i)$ is the same for all $i\in
I$. Let us call it simply $\tau'$.

Let us first consider the case where $\Phi$ is a contraction.

\begin{lem} \label{tdic}
We have a cartesian diagram
\[\comdia{\displaystyle\coprod_{i\in I}\ol{M}(V,\tau',\beta_i)}
{}{\ol{M}(V,\sigma',\beta')} {}{}{}
{\MM(\tau')}{\MM(\Phi')}{\MM(\sigma')}\]
of algebraic $k$-stacks. Moreover,
\[\MM(\Phi')\upsh J(V,\sigma',\beta')=\sum_{i\in
I}J(V,\tau',\beta').\]
\end{lem}
\begin{pf}
The first fact follows immediately from the definitions. The second
fact is \cite{BF} Proposition~7.2.
\end{pf}

Axiom~V will follow by putting Lemma~\ref{tdic} and
Proposition~\ref{evfa} together as follows.
By Lemma~\ref{tdic} all squares in the following diagram are cartesian.
\[\begin{array}{ccccc}
{\displaystyle\coprod_{i\in I}\ol{M}(V,\tau',\beta_i)}&
\stackrel{h}{\longrightarrow} &
\ol{M}(\tau)\times_{\ol{M}(\sigma)}\ol{M}(V,\sigma',\beta') &
\longrightarrow &
\ol{M}(V,\sigma',\beta') \\
\ldiag{} & & \ldiag{} & & \rdiag{a} \\
\MM(\tau') & \stackrel{l}{\longrightarrow} &
\ol{M}(\tau)\times_{\ol{M}(\sigma)}\MM(\sigma') &
\longrightarrow & \MM(\sigma')\\
& \sediag{s} & \ldiag{} & & \rdiag{s} \\
&&\ol{M}(\tau)& \stackrel{\ol{M}(\Phi)}{\longrightarrow} &
\ol{M}(\sigma) \end{array}\]
So we may calculate as follows.
\begin{eqnarray*}
\ol{M}(\Phi)\upsh J(V,\sigma',\beta')  & = &
a\upst s\upst[\ol{M}(\Phi)]\cdot J(V,\sigma',\beta')\\
&=& a\upst l\lst[\MM(\Phi')]\cdot J(V,\sigma',\beta')\\
 \mbox{(by Proposition~\ref{evfa})}& & \\
&=& h\lst\MM(\Phi')\upsh J(V,\sigma',\beta')\\
&=& h\lst\sum_{i\in I} J(V,\tau',\beta_i)
\end{eqnarray*}
by Lemma~\ref{tdic}. This is the context of Axiom~V.

The case that $\Phi$ is of type forgetting a tail is similar. Instead
of Lemma~\ref{tdic} one uses Axiom~IV, and Proposition~\ref{evfa} is
replaced by Proposition~\ref{evfat}.

This finishes the proof of
Axiom~V and hence the proof of Theorem~\ref{ot}.

\subsection{Homogeneous Spaces}

In the case where $V$ is a generalized flag variety, we can give a
more explicit construction of Gromov-Witten invariants as follows.

\subsubsection{Curves and Principal Bundles}

For a smooth algebraic $k$-group $G$ with Lie algebra $\Gg$, we denote
by
\[\HH^1(\tau,G)\]
the $k$-stack of $G$-torsors on $\tau$-marked prestable curves. More
precisely, for a $k$-scheme $T$, the category
$\HH^1(\tau,G)(T)$ is the category of pairs $(C,E)$,
where $C=(C_v)_{v\in V\t}$ is a $\tau$-marked prestable curve over
$T$, giving rise to a morphism $f:T\rightarrow\MM(\tau)$, and $E$ is a
$G$-torsor on $f\upst\cC(\tau)$.

Let $(C,E)$ be such a pair. Denote by $E_v$, for $v\in V\t$, the
$G$-bundle induced by $E$ on $C_v$. We call
\[\deg_v(E)=\deg(E_v)=\deg(E_v\times_{G,Ad}\Gg)\]
the {\em degree }of $E$ at the vertex $v\in V\t$. The degree thus
defines a $\zz_{\geq0}$-structure on $\tau$, which is locally constant
on $T$. (See \cite{BM}, Definition~1.6, for $\zz_{\geq0}$-structures.)

In this way, we get for every $\zz_{\geq0}$-structure $\alpha$ on
$\tau$ an open and closed substack
$\HH^1_{\alpha}(\tau,G)
\subset\HH^1(\tau,G)$,
the substack of $G$-torsors of degree $\alpha$.

\begin{prop} \label{h1sd}
For every $\zz_{\geq0}$-structure $\alpha$ on $\tau$ the stack
$\HH^1_{\alpha}(\tau,G)$ is an algebraic $k$-stack.
The canonical morphism
\[\HH^1_{\alpha}(\tau,G)\longrightarrow\MM(\tau)\]
is smooth of relative dimension
\[-\chi(\tau)\dim G-\alpha(\tau),\]
where $\alpha(\tau)=\sum_{v \in V\t}\alpha(v)$.
\end{prop}
\begin{pf}
To prove that $\HH^1(\tau,G)$ is algebraic, choose a suitable
embedding $G\hookrightarrow GL_n$ to reduce the case of $G$-bundles to the
case of vector bundles, for which it is well-known. The smoothness of
$\HH^1(\tau,G)$ follows from the fact that
$H^2(C,E\times_{G,Ad}\Gg)=0$ for any $G$-torsor $E$ on a $\tau$-marked
prestable curve $C$. The dimension of $\HH^1(\tau,G)$ is equal to
\begin{eqnarray*}
-\chi(E\times_{G,Ad}\Gg)& = &
-\deg(E\times_{G,Ad}\Gg)-\chi(\O_C)\rk(E\times_{G,Ad}\Gg) \\
&=& -\alpha(\tau)-\chi(\tau)\dim G
\end{eqnarray*}
by Riemann-Roch.
\end{pf}

\subsubsection{Maps to $G/P$}

Now let $G$ be a reductive algebraic group over $k$ and $P$ a
parabolic subgroup of $G$. Then $G/P$ is a smooth projective variety
over $k$. Let us assume for simplicity that $G$ is split over $k$. The
morphism $G\rightarrow G/P$ is a principal $P$-bundle, which we shall
denote by $F$.

Let $U_1,\ldots,U_r$ be the elementary representations of $P$ over
$k$, $V_1,\ldots,V_r$ the corresponding vector bundles on $G/P$ and
$L_1,\ldots,L_r$ their determinants. For every $i=1,\ldots,r$ we have
\[V_i=F\times_PU_i.\]
Note that $\Pic(G/P)\otimes\qq$ is spanned by $L_1,\ldots,L_r$ and
that $L_1^{-1}\otimes\ldots\times L_r^{-1}$ is ample.

Let $H_2(G/P)^+$ be the set of homomorphisms of abelian groups
$\psi:\Pic(G/P)\rightarrow\zz$, which are non-negative on ample line
bundles. Then we get a canonical injection
\begin{eqnarray*}
H_2(G/P)^+ & \longrightarrow & (\zz_{\geq0})^r\\
\psi       & \longmapsto     & (\psi(L_1^{-1}),\ldots,\psi(L_r^{-1})).
\end{eqnarray*}
Using this injection we shall think of classes in $H_2(G/P)^+$ as
$r$-tuples of non-negative integers.

Let $\Gg$ and $\Pp$ be the Lie algebras of $G$ and $P$,
respectively. We will consider these only as adjoint representations,
ignoring the Lie algebra structure. Denote by $\Pp$ also the induced
vector bundle
\[F\times_{P,Ad}\Pp\]
on $G/P$. Evaluating on the inverse of its determinant defines a
morphism
\begin{eqnarray*}
\deg:H_2(G/P)^+ & \longrightarrow & \zz_{\geq0} \\
\psi            & \longmapsto     & \psi(\det(\Pp)^{-1}).
\end{eqnarray*}
This morphism has the property that $\deg(\psi)=0$ implies $\psi=0$.

\begin{rmk}
We have $\det\Pp\cong\omega_{G/P}$. In particular,
$\deg\psi=-\psi(\omega_{G/P})$.
\end{rmk}

Now fix an $H_2(G/P)^+$-graph $(\tau,\beta)$, with underlying modular
graph $\tau$. Let $(\tilde{\tau},\tilde{\beta})$ be the $H_2(G/P)^+$-graph
obtained by contracting all edges of $\tau$.

Consider the algebraic $k$-stacks $\HH^1(\tau,G)$ and
$\HH^1(\tau,P)$. Since $G$ is reductive, any $G$-torsor
on a curve has degree zero, and thus
\[\HH^1(\tau,G)\longrightarrow\MM(\tau)\]
is smooth of relative dimension
\[-\chi(\tau)\dim G.\]

If $E$ is a $P$-torsor, then associated to $U_1,\ldots,U_r$ we have
associated vector bundles $E_i=E\times_PU_i$, for $i=1,\ldots,r$, and
thus we may associate to $E$ the {\em multi-degree }
\[\multdeg(E)=(-\deg(E_1),\ldots,-\deg(E_r)).\]
Let $\HH^1_{\beta}(\tau,P)$ be the open and closed
substack of $\HH^1(\tau,P)$ of $P$-torsors whose
multi-degree is equal to $\beta$.

Let $\alpha=\deg\beta$ be the $\zz_{\geq0}$-structure on $\tau$
associated to $\beta$. Then we have
\[\HH^1_{\beta}(\tau,P)
\subset\HH^1_{-\alpha}(\tau,P),\]
so that by Proposition~\ref{h1sd} the stack
$\HH^1_{\beta}(\tau,P)$ is smooth of relative dimension
\[-\chi(\tau)\dim P-\beta(\tau)(\omega_{G/P})\]
over $\MM(\tau)$.

Now let $\MM(G/P,\tau,\beta)$ be the stack of maps from $\tau$-marked
prestable curves to $G/P$ of class $\beta$. More precisely, for a
$k$-scheme $T$, the objects of $\MM(G/P,\tau,\beta)(T)$ are triples
$(C,x,f)$, where $(C,x)$ is a $\tau$-marked prestable curve over $T$
and $f=(f_v)_{v\in V\t}$ is a family of $k$-morphisms
$f_v:C_v\rightarrow G/P$ such that
\begin{enumerate}
\item for all $i\in F\t$ we have
$f_{\del(i)}(x_i)=f_{\del(j\t(i))}(x_{j\t(i)})$,
\item for all $v\in V\t$ we have ${f_v}\lst[C_v]=\beta(v)$.
\end{enumerate}

\begin{rmk}
If $(\tau,\beta)$ is stable, then $\ol{M}(G/P,\tau,\beta)$ is an open
substack of $\MM(G/P,\tau,\beta)$.
\end{rmk}

Note that $G^{V_{\tilde{\tau}}}$ acts on $\MM(G/P,\tau,\beta)$ as
follows. An element $(g_w)_{w\in V_{\tilde{\tau}}}$ of
$G^{V_{\tilde{\tau}}}$ takes $(C,x,(f_v)_{v\in V\t})$ to
$(C,x,(g_{\phi(v)}\comp f_v)_{v\in V\t})$, where
$\phi:\tau\rightarrow\tilde{\tau}$ is the structure contraction.  Let
$$\MM(G/P,\tau,\beta)/G^{V_{\tilde{\tau}}}$$ be the
stack-theoretic quotient of this action. This is an abuse of notation,
since this is a left and not a right action.

We shall let $G^{V_{\tilde{\tau}}}$ act trivially on $\MM(\tau)$
and denote by
\[\MM(\tau)/G^{V_{\tilde{\tau}}}\]
the quotient.

\begin{prop} \label{mcd}
There is a natural cartesian diagram of algebraic $k$-stacks
\[\comdia{\MM(G/P,\tau,\beta)/G^{V_{\tilde{\tau}}}}{\kappa
}{\HH^1_{\beta}(\tau,P)}{\eta}{
}{}{\MM(\tau)/G^{V_{\tilde{\tau}}}}{
\iota}{\HH^1(\tau,G).}\] The vertical maps are
representable, the horizontal maps are local immersions.
\end{prop}
\begin{pf}
This is essentially the fact that a map to $G/P$ is the same as a
principal $P$-bundle with a trivialization of the associated
$G$-bundle.
\end{pf}

The morphism $\iota$ is a local regular immersion with normal bundle
$R^1\pi\lst\O\otimes\Gg$.  Thus the normal cone $C(\tau,\beta)$ of
$\MM(G/P,\tau,\beta)/G^{V_{\tilde{\tau}}}$ in
$\HH^1_{\beta}(\tau,P)$ is a cone in
\[\Nn(\tau,\beta)=\eta\upst R^1\pi\lst\O\otimes\Gg.\]
Pulling back to $\MM(G/P,\tau,\beta)$ and, if $(\tau,\beta)$ is
stable, to $\ol{M}(G/P,\tau,\beta)$ defines
$G^{V_{\tilde{\tau}}}$-equivariant cones, which we shall still
denote $C(\tau,\beta)$, inside equivariant vector bundles, which we
shall still denote by $\Nn(\tau,\beta)$.

Let us now assume that $(\tau,\beta)$ is stable. Then we may intersect
the cone $C(\tau,\beta)$ over $\ol{M}(G/P,\tau,\beta)$ with the zero
section of the vector bundle $\Nn(\tau,\beta)$, to define a cycle
class
\[J(\tau,\beta)\in A_{\dim(G/P,\tau,\beta)}(\ol{M}(G/P,\tau,\beta))\]
with rational coefficients. Note that $C(\tau,\beta)$ is pure of the
correct dimension, since it is constructed as a normal cone inside a
smooth stack of the correct dimension.

\begin{prop} \label{tmt}
The collection of cycle classes $J(\tau,\beta)$ is the orientation of
$\ol{M}$ over $\GG_s(G/P)$ defined using the intrinsic normal cone.
\end{prop}
\begin{pf}
This follows from \cite{BF} Example~7.6, since
\[\dual{(R\pi\lst f\upst T_{G/P})}=\kappa\upst
L\com_{\HH^1_{\beta}(\tau,P)/\HH^1(\tau,G)}.\]
\end{pf}

\begin{rmk}
As a corollary we get that the orientation classes $J(\tau,\beta)$ are
$G^{V_{\tilde{\tau}}}$-invariant. The same is then true for the
Gromov-Witten invariants.
\end{rmk}


\begin{thebibliography}{1}

\bibitem{BF}
K.~Behrend and B.~Fantechi.
\newblock The intrinsic normal cone.
\newblock Preprint, 1996.

\bibitem{BM}
K.~Behrend and Yu. Manin.
\newblock Stacks of stable maps and {G}romov-{W}itten invariants.
\newblock To appear in {\em Duke Mathematical Journal}.

\bibitem{fgaIV}
A.~Grothendieck.
\newblock Techniques de construction et th\'eor\`emes d'existence en
  g\'eom\'etrie alg\'ebrique {IV}: Les sch\'emas de {H}ilbert.
\newblock {\em S\'eminaire Bourbaki}, 13e ann\'ee(221), 1960--61.

\bibitem{knudsen}
F.~Knudsen.
\newblock The projectivity of the moduli space of stable curves, {I}{I}: The
  stacks ${M}_{g,n}$.
\newblock {\em Math.\ Scand.}, 52:161--199, 1983.

\bibitem{KM}
M.~Kontsevich and Yu. Manin.
\newblock {G}romov-{W}itten classes, quantum cohomology, and enumerative
  geometry.
\newblock {\em Commuications in Mathematical Physics}, 164:525--562, 1994.

\end{thebibliography}

\begin{flushleft}
{\tt behrend@math.ubc.ca}
\end{flushleft}

\end{document}